\documentclass[12pt]{article}
\usepackage{amssymb,fullpage,enumerate,amsfonts,amsmath,latexsym}

\newtheorem{definition}{Definition}[section]
\newtheorem{theorem}[definition]{Theorem}
\newtheorem{lemma}[definition]{Lemma}
\newtheorem{corollary}[definition]{Corollary}
\newtheorem{remark}[definition]{Remark}

\newtheorem{note}[definition]{Note}

\newtheorem{assumption}[definition]{Assumption}

\def\K{\mathbb K}

\def\Z{\mathbb Z}
\def\K{\mathbb K}

\newcommand{\uqhat}{U_q(\widehat{\mathfrak{sl}}_2)}

\begin{document}

\title{ \bf Raising/Lowering Maps and Modules for the Quantum Affine Algebra $\uqhat$}
\author {Darren Funk-Neubauer}
\date{}
\maketitle 

\begin{abstract}
Let $V$ denote a finite dimensional vector space over an algebraically closed field.  Let $U_0, U_1,\ldots, U_d$ denote a sequence of nonzero subspaces whose direct sum is $V$.  Let $R:V\to V$ and $L:V\to V$ denote linear transformations with the following properties:  for $0 \leq i \leq d$, $R U_i \subseteq U_{i+1}$ and $L U_i \subseteq U_{i-1}$ where $U_{-1}=0$, $U_{d+1}=0$;  for $ 0 \leq i \leq d/2$, the restrictions $R^{d-2i}|_{U_i}: U_i\to U_{d-i}$ and $L^{d-2i}|_{U_{d-i}}: U_{d-i}\to U_i$ are bijections;  the maps $R$ and $L$ satisfy the cubic $q$-Serre relations where $q$ is nonzero and not a root of unity.  Let $K:V\to V$ denote the linear transformation such that $(K-q^{2i-d}I)U_i=0$ for $0 \leq i \leq d$.  We show that there exists a unique $\uqhat$-module structure on $V$ such that each of $R-e^{-}_1$, $L-e^{-}_0$, $K-K_0$, and $K^{-1}-K_1$ vanish on $V$, where $e^{-}_1, e^{-}_0, K_0, K_1$ are Chevalley generators for $\uqhat$.  We determine which $\uqhat$-modules arise from our construction.  
\end{abstract}

\section{Introduction}

\noindent
A quantum affine algebra is a $q$-analogue of the universal enveloping algebra of a Kac-Moody Lie algebra of affine type.  Quantum affine algebras were first introduced and studied by V.G. Drinfeld \cite{Drinfeld1, Drinfeld2} and M. Jimbo \cite{Jimbo1, Jimbo2} in relation to the Yang-Baxter equation of mathematical physics.  Since then quantum affine algebras have played an important role in various areas of mathematics and physics, for example see \cite{Beck}, \cite{Kang}, \cite{Mansour}, \cite{Misra}.  In this paper we will be concerned with the quantum affine algebra $\uqhat$.  \\

\noindent
We study the finite dimensional modules for $\uqhat$.  In \cite{cp} V. Chari and A. Pressley classified the finite dimensional irreducible $\uqhat$-modules up to isomorphism.  These modules were further studied in \cite{cp2}, \cite{Thoren}.  However, a complete classification of all finite dimensional $\uqhat$-modules is still unknown.  In the present paper we obtain a result that may help in this classification.  We now introduce some notation and recall the definition of $\uqhat$. \\

\noindent 
Throughout this paper $\K$ will denote an algebraically closed field.  We fix a nonzero scalar $q\in \K$ that is not a root of unity.  We will use the following notation:
\begin{eqnarray*}
\lbrack n \rbrack = {{q^n-q^{-n}} \over {q-q^{-1}}},
\qquad \qquad n=0,1,\ldots.
\label{eq:brackndef}
\end{eqnarray*}

\begin{definition}
\label{def:uq}
\rm
\cite[Definition 2.2]{cp}
The quantum affine algebra $\uqhat$ is the unital associative $\K$-algebra with generators $e^{\pm}_i$, $K_i^{{\pm}1}$, $i\in \lbrace 0,1\rbrace$, which satisfy the following relations:  
\begin{eqnarray}
K_iK^{-1}_i&=&
K^{-1}_iK_i =  1,
\label{eq:buq1}
\\
K_0K_1&=& K_1K_0,
\label{eq:buq2}
\\
K_ie^{\pm}_iK^{-1}_i &=& q^{{\pm}2}e^{\pm}_i,
\label{eq:buq3}
\\
K_ie^{\pm}_jK^{-1}_i &=& q^{{\mp}2}e^{\pm}_j, \qquad i\not=j,
\label{eq:buq4}
\\
e^{+}_ie^{-}_i-e^{-}_ie^{+}_i  &=& {{K_i-K^{-1}_i}\over {q-q^{-1}}}, 
\label {eq:buq5}
\\
e^{\pm}_0e^{\mp}_1 &=& e^{\mp}_1e^{\pm}_0,
\label{eq:buq6}
\end{eqnarray} 
\begin{eqnarray}
(e^{\pm}_i)^3e^{\pm}_j - \lbrack 3 \rbrack (e^{\pm}_i)^2e^{\pm}_je^{\pm}_i + \lbrack 3 \rbrack e^{\pm}_ie^{\pm}_j(e^{\pm}_i)^2 - e^{\pm}_j (e^{\pm}_i)^3 =0, \qquad i\not=j.
\label{eq:buq7}
\end{eqnarray}
We call $e^{\pm}_i$, $K_i^{{\pm}1}$ the {\it Chevalley generators} for $\uqhat$ and refer to (\ref{eq:buq7}) as the {\it cubic $q$-Serre relations}.
\end{definition}

\noindent
We now give two definitions, state our main result, and then make some comments concerning its significance.  
  
\begin{definition}
\label{def:decomp}
\rm
Let $V$ denote a vector space over $\K$ with finite positive dimension.  By a \\
{\it decomposition} of $V$ we mean a sequence $U_0, U_1, \ldots, U_d$ consisting of nonzero subspaces of $V$ such that $V=\sum_{i=0}^{d} U_i$ (direct sum).  For notational convenience we set $U_{-1}:=0, U_{d+1}:=0$.
\end{definition}

\begin{definition}
\label{def:K}
\rm
Let $V$ denote a vector space over $\K$ with finite positive dimension.  Let $U_0, U_1, \ldots, U_d$ denote a decomposition of $V$.  Let $K:V\to V$ denote the linear transformation such that, for $0 \leq i \leq d$, $U_i$ is an eigenspace for $K$ with eigenvalue $q^{2i-d}$.  We refer to $K$ as {\it the linear transformation that corresponds to the decomposition $U_0, U_1, \ldots, U_d$}.  
\end{definition}

\begin{note}
\label{def:K^-1}
\rm
\noindent
With reference to Definition \ref{def:K}, we note that $K$ is invertible.  Moreover, for $0 \leq i \leq d$, $U_{i}$ is the eigenspace for $K^{-1}$ with eigenvalue $q^{d-2i}$.  We observe that $K^{-1}$ is the linear transformation that corresponds to the decomposition  $U_d, U_{d-1}, \ldots, U_0$.
\end{note}

\noindent
We will be concerned with the following situation.

\begin{assumption}
\label{ass}
\rm
Let $V$ denote a vector space over $\K$ with finite positive dimension.  Let $U_0, U_1, \ldots, U_d$ denote a decomposition of $V$.  Let $K$ denote the linear transformation that corresponds to $U_0, U_1, \ldots, U_d$ as in Definition \ref{def:K}.   Let $R:V\to V$ and $L:V\to V$ be linear transformations such that 
\begin{enumerate}
\item $R U_i \subseteq U_{i+1} \qquad (0 \leq i \leq d-1), \qquad  R U_d=0$,
\item $L U_i \subseteq U_{i-1} \qquad (1 \leq i \leq d), \qquad  L U_0=0$,
\item for $ 0 \leq i \leq d/2$ the restriction $R^{d-2i}|_{U_i}: U_i\to U_{d-i}$ is a bijection,
\item for $0 \leq i \leq d/2$ the restriction $L^{d-2i}|_{U_{d-i}}: U_{d-i}\to U_i$ is a bijection, 
\item $R^3L-\lbrack 3 \rbrack R^2LR+\lbrack 3 \rbrack RLR^2-LR^3=0$,
\item $L^3R-\lbrack 3 \rbrack L^2RL+\lbrack 3 \rbrack LRL^2-RL^3=0$.
\end{enumerate}  
\end{assumption}

\noindent
We now state our main result.

\begin{theorem}
\label{thm:main}
Adopt Assumption \ref{ass}.  Then there exists a unique $\uqhat$-module structure on $V$ such that $(R-e^{-}_1)V=0$, $(L-e^{-}_0)V=0$, $(K-K_0)V=0$, $(K^{-1}-K_1)V=0$, where $e^{-}_1$, $e^{-}_0$, $K_0$, $K_1$ are Chevalley generators for $\uqhat$ as in Definition \ref{def:uq}. 
\end{theorem}

\noindent
The proof of Theorem \ref{thm:main} will take up most of the paper until Section 9.  In Sections 10 and 11 we determine which $\uqhat$-modules arise from the construction of Theorem \ref{thm:main}.  

\begin{remark}
\rm
Not all finite dimensional $\uqhat$-modules arise from the construction of Theorem \ref{thm:main}.  However as we will see, every finite dimensional $\uqhat$-module is a direct sum of submodules, each of which arises from Theorem \ref{thm:main} up to a routine normalization.  Thus Theorem \ref{thm:main} can be viewed as a step towards the classification of the finite dimensional $\uqhat$-modules.  
\end{remark}

\begin{remark}
\rm
The proof of Theorem \ref{thm:main} involves modifying a construction used in \cite{Benkart} and \cite{Ito.Ter}.  The construction originally arose from the study of tridiagonal pairs.  According to \cite[Definition 1.1]{Ito} a tridiagonal pair is an ordered pair $(A,A^*)$ of diagonalizable linear transformatons on a finite dimensional vector space $V$ such that (i) the eigenspaces of $A$ (resp. $A^*$) can be ordered as $V_0, V_1, \ldots, V_d$ (resp. $V^*_0, V^*_1, \ldots, V^*_d$) with $A^*V_i \subseteq V_{i-1} + V_{i} + V_{i+1}$ (resp. $AV^*_i \subseteq V^*_{i-1} + V^*_{i} + V^*_{i+1}$) for $0 \leq i \leq d$;  (ii) there are no nonzero proper subspaces of $V$ which are invariant under both $A$ and $A^*$.  See \cite{Curtin}, \cite{shape}, \cite{Ito.Ter} for connections between tridiagonal pairs and $\uqhat$. 
\end{remark}

\begin{remark}
\rm
In \cite{Benkart} G. Benkart and P. Terwilliger determine the finite dimensional irreducible modules for the standard Borel subalgebra of $\uqhat$.  The authors adopt Assumption \ref{ass}(i),(ii),(v),(vi) but replace Assumption \ref{ass}(iii),(iv) with the assumption that $V$ is irreducible as a $(K,R,L)$-module.  From this assumption they obtain a $\uqhat$-module structure on $V$ as in Theorem \ref{thm:main}.  The $\uqhat$-module structure that they obtain is irreducible while the $\uqhat$-module structure given by our Theorem \ref{thm:main} is not necessarily irreducible.  As far as we know Theorem \ref{thm:main} does not imply the result in \cite{Benkart} nor does the result in \cite{Benkart} imply Theorem \ref{thm:main}.  
\end{remark}

\begin{remark}
\rm
In the proof of Theorem \ref{thm:main} we will need a number of lemmas that are similiar to lemmas appearing in \cite{Benkart}.  The assumptions in Theorem \ref{thm:main} and in \cite{Benkart} are similiar except that in \cite{Benkart} the authors assume $V$ is irreducible as a $(K,R,L)$-module and we are not making this assumption.  However,  in \cite{Benkart} this irreducibility condition is not used in the proofs of the lemmas we require.  In such cases we simply cite \cite{Benkart} without proof.  
\end{remark}

\section{Preliminaries}
\noindent
In this section we make a few observations about Assumption \ref{ass}.

\begin{note}
\label{note:assumptions} 
\rm
\noindent
With reference to Assumption \ref{ass}, if we replace $U_i$ by $U_{d-i}$ for $0 \leq i \leq d$, and replace $(K,R,L)$ by $(K^{-1},L,R)$, then the assumption is still satisfied.   
\end{note}

\begin{lemma}
\label{thm:dim}
With reference to Assumption \ref{ass}, for $0 \leq i \leq  d/2$ and $0 \leq j \leq d-2i$, the restriction $R^{j}|_{U_i}:U_i\to U_{i+j}$ is an injection.
\end{lemma}

\noindent
{\it Proof:} Immediate from Assumption \ref{ass}(iii).   
\hfill $\Box $ \\

\begin{lemma}
\label{thm:KLR}
With reference to Assumption \ref{ass}, the following (i),(ii) hold.
\begin{enumerate}
\item $KR=q^2RK$,
\item $KL=q^{-2}LK$.
\end{enumerate}
\end{lemma}

\noindent
{\it Proof:} Immediate from Assumption \ref{ass} (i),(ii) and Definition \ref{def:K}.
\hfill $\Box$ \\

\section{An outline of the proof for Theorem \ref{thm:main}}

\noindent
Our proof of Theorem \ref{thm:main} will consume most of the paper from Section
 4 to Section 9.  Here we sketch an overview of the argument.

\medskip
\noindent
We begin by adopting Assumption \ref{ass}.  To start the construction of 
the $\uqhat$-action on $V$ we require that the linear transformations $R-e^-_1$, $L-e^-_0$, $K^{\pm 1}-K^{\pm 1}_0$, and $K^{\pm 1}- K^{\mp 1}_1$ vanish on $V$.  This gives the actions of the elements $e^-_1, e^-_0, K^{\pm 1}_0, K^{\pm 1}_1$ on $V$.  We define the actions of $e^+_0, e^+_1$ on $V$ as follows.  First we prove that $K+R$ is diagonalizable on $V$.  Then we show that the set of distinct eigenvalues of $K+R$ on $V$ is $\{\, q^{2i-d} \, | \, 0 \leq i \leq d \, \}$.  For $0 \leq i \leq d$ we let $V_i$ denote the eigenspace of $K+R$ on $V$ associated with the eigenvalue $q^{2i-d}$.  So $V_0, V_1,\ldots, V_d$  is a decomposition of $V$.  Next we define the subspaces $W_i$ as follows.
\begin{eqnarray*}
W_i =(U_0 + \cdots + U_i) \cap (V_0+\cdots + V_{d-i})  \qquad (0 \leq i \leq d).
\end{eqnarray*}
We show that $W_0, W_1, \ldots, W_d$ is a decomposition of $V$.  Then we apply Note \ref{note:assumptions} to the above argument to obtain the following results.  $K^{-1}+L$ is diagonalizable on $V$.  The set of distinct eigenvalues of $K^{-1}+L$ on $V$ is $\{\, q^{d-2i} \, | \, 0 \leq i \leq d \, \}$.  For $0 \leq i \leq d$ we let $V^*_i$ denote the eigenspace of $K^{-1}+L $ on $V$ associated with the eigenvalue $q^{d-2i}$.  So $V^*_0, V^*_1, \ldots, V^*_d$  is a decomposition of $V$.  Next we define the subspaces $W^*_i$ as follows.
\begin{eqnarray*}
W^*_i =(U_{i} + \cdots + U_d) \cap (V^*_{d-i}+\cdots + V^*_{d})  \qquad (0 \leq i \leq d).
\end{eqnarray*}
Then $W^*_0, W^*_1, \ldots, W^*_d$ is a decomposition of $V$.  Next we define the linear transformation $B:V\to V$ (resp. $B^*:V\to V$) such that for $0 \leq i \leq d$, $W_i$ (resp. $W^*_i$) is an eigenspace for $B$ (resp. $B^*$) with eigenvalue $q^{2i-d}$ (resp. $q^{d-2i}$).  We let $e^+_1$ act on $V$ as $I-K^{-1}B$ times $q^{-1}(q-q^{-1})^{-2}$.  We let $e^+_0$ act on $V$ as $I-KB^*$ times $q^{-1}(q-q^{-1})^{-2}$.  We display some relations that are satisfied by $B,B^*,L,R,K,K^{-1}$.  Using these relations we argue that the above actions of $e_0^{\pm}, e_1^{\pm}, K_0^{\pm 1}, K_1^{\pm 1}$ satisfy the defining relations for $\uqhat$.  In this way we obtain the required action of $\uqhat$ on $V$. 

\section{The linear transformation $A$}
\noindent
In this section we define and discuss a linear transformation that will be useful.

\begin{definition}
\label{def:A}
\rm
With reference to Assumption \ref{ass}, let $A:V\to V$ denote the following linear transformation:
\begin{eqnarray*} 
 A=K+R.
\end{eqnarray*}
\end{definition}

\begin{lemma}
\label{thm:AKR}
With reference to Definition \ref{def:A} and Assumption \ref{ass}, the following (i),(ii) hold.
\begin{enumerate}
\item For $0 \leq i \leq d$ the action of $A-q^{2i-d}I$ on $U_i$ coincides with the action of $R$ on $U_i$,
\item $(A-q^{2i-d}I)U_i \subseteq U_{i+1}$, \qquad $0 \leq i \leq d$.
\end{enumerate}
\end{lemma}

\noindent
{\it Proof:}  Immediate from Definition \ref{def:A}, Definition \ref{def:K}, and Assumption \ref{ass}(i).  \hfill $\Box$ \\

\begin{lemma}
\label{thm:Adiag}
{\rm \cite[Lemma 4.13]{Benkart}}
With reference to Definition \ref{def:A} and Assumption \ref{ass}, the following holds. 
$A$ is diagonalizable on $V$ and the set of distinct eigenvalues of $A$ is \\ $\{\,q^{2i-d} \, |\, 0 \leq i \leq d \, \}$.  Moreover, for $0 \leq i \leq d$, the dimension of the eigenspace for $A$ associated with $q^{2i-d}$ is equal to the dimension of $U_i$. 
\end{lemma}

\begin{definition}
\label{def:V}
\rm
With reference to Definition \ref{def:A} and Lemma \ref{thm:Adiag}, for $0 \leq i \leq d$ we let $V_i$ denote the eigenspace for $A$ with eigenvalue $q^{2i-d}$.  For notational convenience we set $V_{-1}:=0,V_{d+1}:=0$.  We observe that $V_0, V_1, \ldots, V_d$ is a decomposition of $V$.  
\end{definition}

\begin{lemma}
\label{thm:samedim}
Let the decomposition $U_0, U_1, \ldots, U_d$ be as in Assumption \ref{ass} and let the decomposition $V_0, V_1, \ldots, V_d$ be as in Definition \ref{def:V}.  Then for $0 \leq i \leq d$ the spaces $U_i$, $U_{d-i}$, $V_i$, $V_{d-i}$ all have the same dimension. 
\end{lemma}  

\noindent
{\it Proof:} Immediate from Lemma \ref{thm:Adiag} and Assumption \ref{ass} (iii).
\hfill $\Box$ \\

\begin{definition}
\label{def:rho}
\rm
With reference to Lemma \ref{thm:samedim}, for $0 \leq i \leq d$ we let $\rho_i$ denote the common dimension of $U_i$, $U_{d-i}$, $V_i$, $V_{d-i}$. 
\end{definition} 

\begin{lemma}
\label{thm:rhofacts}
With reference to Definition \ref{def:rho}, the following (i)--(iii) hold.
\begin{enumerate}
\item $\rho_i \neq 0$, \qquad $0 \leq i \leq d$,
\item $\dim(V)=\sum_{i=0}^{d}\rho_i$,
\item $\rho_i=\rho_{d-i}$, \qquad $0 \leq i \leq d$.
\end{enumerate}
\end{lemma}

\noindent 
{\it Proof:}  Immediate by Definition \ref{def:rho} and since $U_0, U_1, \ldots, U_d$ is a decomposition of $V$.
\hfill $\Box$ \\

\begin{lemma} 
\label{thm:UV}
{\rm \cite[Lemma 5.2]{Benkart}}
With reference to Assumption \ref{ass} and Definition \ref{def:V},  
\begin{eqnarray*}
U_i+\cdots+U_d=V_i+\cdots+V_d,  \qquad 0 \leq i \leq d.
\end{eqnarray*}
\end{lemma}

\begin{lemma}
\label{thm:AK^{-1}relation}
{\rm \cite[Lemma 5.3]{Benkart}}
With reference to Assumption \ref{ass} and Definition \ref{def:V}, 
\begin{eqnarray}
(K^{-1}-q^{d-2i})V_i \subseteq V_{i+1}, \qquad 0 \leq i \leq d. 
\end{eqnarray}
\end{lemma}

\section{The subspaces $W_i$}

\begin{definition}
\label{def:W}
\rm  
With reference to Assumption \ref{ass} and Definition \ref{def:V}, we define
\begin{eqnarray*}
W_i=(U_0+ \cdots +U_i) \cap (V_0+ \cdots +V_{d-i}),  \qquad 0 \leq i \leq d.
\end{eqnarray*}
For notational convenience we set $W_{-1}:=0, W_{d+1}:=0$.
\end{definition}

\noindent
The goal of this section is to prove the following theorem.

\begin{theorem}
\label{thm:Wfinal}
With reference to Definition \ref{def:W}, the sequence $W_0, W_1, \ldots, W_d$ is a decomposition of $V$. 
\end{theorem}

\noindent
We prove Theorem \ref{thm:Wfinal} in three steps.  First, we show the sum $\sum_{i=0}^dW_i$ is direct.  Second, we show $V=\sum_{i=0}^dW_i$.  Finally, we show $W_i \neq 0$ for $0 \leq i \leq d$. \\

\noindent
The following definition and the next few lemmas will be useful in proving the sum $\sum_{i=0}^dW_i$ is direct.

\begin{definition}
\label{def:Wij}
\rm
With reference to Assumption \ref{ass} and Definition \ref{def:V}, we define
\begin{eqnarray*}
W(i,j)=\biggl(\sum_{h=0}^{i}U_h\biggr) \, \cap \, \biggl(\sum_{h=0}^{j}V_h\biggr), \qquad -1 \leq i, j \leq d+1.
\end{eqnarray*}
With reference to Definition \ref{def:W}, note that $W(i,d-i)=W_i$ for $0 \leq i \leq d$.
\end{definition}

\begin{lemma}
\label{thm:AK^{-1}Wij}
With reference to Assumption \ref{ass}, Definition \ref{def:A}, and Definition \ref{def:Wij}, the following (i),(ii) hold.
\begin{enumerate}
\item $(A-q^{2j-d}I)W(i,j) \subseteq W(i+1,j-1)$, \qquad $0 \leq i,j \leq d$,
\item $(K^{-1}-q^{d-2i}I)W(i,j) \subseteq W(i-1,j+1)$, \qquad $0 \leq i,j \leq d$.
\end{enumerate}
\end{lemma}

\noindent
{\it Proof:} (i) Using Definition \ref{def:Wij} and  Definition \ref{def:V}, we have $(A-q^{2j-d}I)W(i,j) \subseteq \sum_{h=0}^{j-1}V_h$.  Using Definition \ref{def:Wij} and Lemma \ref{thm:AKR}(ii) we have $(A-q^{2j-d}I)W(i,j) \subseteq \sum_{h=0}^{i+1}U_h$.  Combining these facts we obtain the desired result. \\
(ii) Using Definition \ref{def:Wij} and Lemma \ref{thm:AK^{-1}relation}, we have $(K^{-1}-q^{d-2i}I)W(i,j) \subseteq \sum_{h=0}^{j+1}V_h$.  Using Definition \ref{def:Wij} and Note \ref{def:K^-1} we have $(K^{-1}-q^{d-2i}I)W(i,j) \subseteq \sum_{h=0}^{i-1}U_h$.  Combining these facts we obtain the desired result.
\hfill $\Box$ \\

\begin{lemma}
\label{thm:wij=0}
With reference to Definition \ref{def:Wij},
\begin{eqnarray*}
W(i,d-1-i)=0,  \qquad \qquad 0 \leq i \leq d-1.
\end{eqnarray*}
\end{lemma}

\noindent
{\it Proof:} Define $T=\sum_{i=0}^{d-1}W(i,d-1-i)$.  It suffices to show $T=0$.  By Lemma \ref{thm:AK^{-1}Wij}(ii) we find $K^{-1}T \subseteq T$.  Recall that $K^{-1}$ is diagonalizable on $V$ and so $K^{-1}$ is diagonalizable on $T$.  Also, $U_j \cap T$ $(0 \leq j \leq d)$ are the eigenspaces for $K^{-1}|_{T}$.  Thus $T=\sum_{j=0}^{d}(U_j \cap T)$ (direct sum).  Suppose, towards a contradiction, $T \neq 0$.  Then there exists $j$ $(0 \leq j \leq d)$ such that $U_j \cap T \neq 0$.  Define $t:=\mbox{min}\{\,j \, |\, 0 \leq j \leq d\, , \, U_j \cap T \neq 0 \,\}$ and $r:=\mbox{max}\{\,j \, |\, 0 \leq j \leq d\, , \, U_j \cap T \neq 0 \,\}$.  Of course $t \leq r$.  We will now show 
\begin{eqnarray}
\label{numb1}
r+t \geq d.
\end{eqnarray}
If $d/2  <  t$ then (\ref{numb1}) holds since $t \leq r$.  So now assume $0 \leq t \leq d/2$.  Let $x \in U_t \cap T$ such that $x \neq 0$.  By Assumption \ref{ass}(i), we find $R^{d-2t}x \in U_{d-t}$.  Also, by Assumption \ref{ass}(iii), we find $R^{d-2t}x \neq 0$.  By Lemma \ref{thm:AK^{-1}Wij}(i), we find $AT \subseteq T$.  Using this and Lemma \ref{thm:AKR}(i), we find $R^{d-2t}x \in T$.  So $R^{d-2t}x \in U_{d-t} \cap T$.  Combining these facts we find $U_{d-t} \cap T \neq 0$.  This shows (\ref{numb1}).  Define $y:=\mbox{max}\{\,j \, |\, 0 \leq j \leq d-1\, , \, W(j,d-1-j) \neq 0 \,\}$.  By the definition of $T$ we find $T \subseteq U_0+ \cdots +U_y$ and so
\begin{eqnarray}
\label{numb2}
y \geq r.
\end{eqnarray}
We will now show 
\begin{eqnarray}
\label{numb3}
d-y \geq t+1.
\end{eqnarray}
By the definition of $y$ we have $W(y,d-1-y) \neq 0$.  By Definition \ref{def:Wij} we have $W(y,d-1-y) \subseteq V_0+ \cdots +V_{d-1-y}$.  Using these facts and since $V_0, V_1, \ldots, V_d$ is a decomposition of $V$ we find $W(y,d-1-y) \nsubseteq V_{d-y}+ \cdots +V_d$.  Therefore $T \nsubseteq V_{d-y}+ \cdots +V_d$.  Using this and Lemma \ref{thm:UV} we find $T \nsubseteq U_{d-y}+ \cdots +U_d$ and (\ref{numb3}) follows.  Adding (\ref{numb1}), (\ref{numb2}), and (\ref{numb3}) we find $0 \geq 1$ for a contradiction.  The result follows.
\hfill $\Box$. \\       

\begin{lemma}
\label{thm:directsum}
With reference to Definition \ref{def:W}, the sum $\sum_{i=0}^{d}W_i$ is direct.
\end{lemma}

\noindent
{\it Proof:}  It suffices to show $(W_0+ \cdots +W_{i-1}) \cap W_i=0$ for $1 \leq i \leq d$.  Let $i$ be given.  By Definition \ref{def:W}, $W_0+ \cdots +W_{i-1} \subseteq U_0+ \cdots +U_{i-1}$.  Also by Definition \ref{def:W}, $W_i \subseteq V_0+ \cdots +V_{d-i}$.  By this and Definition \ref{def:Wij} we find $(W_0+ \cdots +W_{i-1}) \cap W_i \subseteq W(i-1,d-i)$.  But $W(i-1,d-i)=0$ by Lemma \ref{thm:wij=0} and so $(W_0+ \cdots +W_{i-1}) \cap W_i=0$.
\hfill $\Box$ \\

\noindent
The following definition and the next few lemmas will be useful in proving $V=\sum_{i=0}^{d}W_i$.

\begin{definition}
\label{def:H_i}
\rm
With reference to Assumption \ref{ass}, we define
\begin{eqnarray*}
H_i=\{\,v \in U_i\, |\, R^{d-2i+1}v=0\,\}, \qquad 0 \leq i \leq d/2.
\end{eqnarray*}
\end{definition}

\begin{lemma}
\label{thm:HVU}
With reference to Assumption \ref{ass}, Definition \ref{def:V}, and Definition \ref{def:H_i},\begin{eqnarray*}
H_i=(V_i+ \cdots +V_{d-i}) \cap U_i,   \qquad 0 \leq i \leq d/2.
\end{eqnarray*}
\end{lemma}

\noindent
{\it Proof:} Immediate from Definition \ref{def:V}, Definition \ref{def:H_i}, and Lemma \ref{thm:AKR}(i). 
\hfill $\Box$ \\

\begin{lemma}
\label{thm:refineU_i}
With reference to Assumption \ref{ass} and Definition \ref{def:H_i}, 
\begin{eqnarray*}
U_i=\sum_{j=0}^{min(i,d-i)}R^{i-j}H_j \qquad (direct \, \, sum), \qquad 0 \leq i \leq d.
\end{eqnarray*}
\end{lemma}

\noindent
{\it Proof:}  Case 1: $0 \leq i \leq d/2$.  The proof is by induction on $i$.  Observe the result holds for $i=0$ since $U_0=H_0$ by Definition \ref{def:H_i} and Assumption \ref{ass}(i).  Next assume $i \geq 1$. By induction and Lemma \ref{thm:dim} we find
\begin{eqnarray}
\label{refineU_i1}
RU_{i-1}=\sum_{j=0}^{i-1}R^{i-j}H_j \qquad  (direct \, \, sum).
\end{eqnarray}
We now show 
\begin{eqnarray}
\label{refineU_i2}
U_i=RU_{i-1}+H_i \qquad (direct \, \, sum).
\end{eqnarray}
Using Assumption \ref{ass}(i) and Definition \ref{def:H_i}, we have $RU_{i-1}+H_i \subseteq U_i$.  We now show $U_i \subseteq RU_{i-1}+H_i$.  Let $x \in U_i$.  By Assumption \ref{ass}(i),(iii) there exists $y \in U_{i-1}$ such that $R^{d-2i+2}y=R^{d-2i+1}x$.  Using this we find $x-Ry \in H_i$.  So $x \in RU_{i-1}+H_i$.  We have now shown equality in (\ref{refineU_i2}).  It remains to show that the sum in (\ref{refineU_i2}) is  direct.  To do this we show $RU_{i-1} \cap H_i=0$.  Let $x \in RU_{i-1} \cap H_i$.  By Definition \ref{def:H_i} we have $R^{d-2i+1}x=0$.  Also, there exists $y \in U_{i-1}$ such that $x=Ry$.  Combining these facts with Assumption \ref{ass}(iii) we find $y=0$ and then $x=0$.  We have now shown the sum in (\ref{refineU_i2}) is direct and this completes the proof of (\ref{refineU_i2}).  Combining (\ref{refineU_i1}) and (\ref{refineU_i2}) we find 
\begin{eqnarray*}
U_i=\sum_{j=0}^{i}R^{i-j}H_j \qquad (direct \, \, sum).
\end{eqnarray*} 
Case 2: $d/2 < i \leq d$.  This case follows immediately from Case 1 and Assumption \ref{ass}(iii).
\hfill $\Box$ \\

\noindent
Recall that $\mbox{End}(V)$ is the $\K$-algebra consisting of all linear transformations from $V$ to $V$.

\begin{definition}
\label{def:scriptD}
\rm
With reference to Definition \ref{def:A}, let $\mathcal{D}$ denote the $\K$-subalgebra of $\mbox{End}(V)$ generated by $A$.  
\end{definition}

\noindent
We will be concerned with the following subspace of $V$.  With reference to Definition \ref{def:H_i} and Definition \ref{def:scriptD}, for $0 \leq i \leq d/2$ we define $\mathcal{D}H_i=\mathrm{span}\{\, Xh\, |\,X\in \mathcal{D}, \,h\in H_i\,\}$.

\begin{lemma}
\label{thm:scriptDR^iH_i}
With reference to Assumption \ref{ass}, Definition \ref{def:H_i}, and  Definition \ref{def:scriptD}, 
\begin{eqnarray}
\mathcal{D}H_i=\sum_{j=0}^{d-2i}R^jH_i \qquad (direct \,\, sum), \qquad  0 \leq i \leq d/2.
\label{scriptDR^iH_i1}
\end{eqnarray}
\end{lemma}

\noindent 
{\it Proof:} Let $i$ be given.  Define $\Delta=\sum_{j=0}^{d-2i}R^jH_i$.  We first show $\mathcal{D}H_i=\Delta$.  Recall $H_i \subseteq U_i$ by Definition \ref{def:H_i}.  By this and Lemma \ref{thm:AKR}(i) we find $\Delta \subseteq \mathcal{D}H_i$.  We now show $\mathcal{D}H_i \subseteq \Delta$.  Since $\mathcal{D}$ is generated by $A$ and since $\Delta$ contians $H_i$ it suffices to show that $\Delta$ is $A$-invariant.  We now show $\Delta$ is $A$-invariant.  For $0 \leq j \leq d-2i$ and $h \in H_i$ we show $AR^jh \in \Delta$.  Using Assumption \ref{ass}(i) and Lemma \ref{thm:AKR}(i) we find $AR^jh \in R^jH_i+R^{j+1}H_i$.  Recall $R^{d-2i+1}H_i=0$ by Definition \ref{def:H_i}.  By these comments we find $AR^jh \in \Delta$.  This completes the proof that $\Delta$ is $A$-invariant and it follows $\mathcal{D}H_i \subseteq \Delta$.  We have now shown $\mathcal{D}H_i=\Delta$.  It remains to show that the sum $\sum_{j=0}^{d-2i}R^jH_i$ is direct.  This follows since $U_0, U_1, \ldots, U_d$ is a decomposition of $V$ and since $R^jH_i \subseteq U_{i+j}$ $(0 \leq j \leq d-2i)$ by Assumption \ref{ass}(i).  
\hfill $\Box$ \\
   
\begin{lemma}
\label{thm:VscriptDH_i}
With reference to Assumption \ref{ass} and Definition \ref{def:scriptD},  
\begin{eqnarray*}
V=\sum_{i=0}^{d/2}\mathcal{D}H_i \qquad (direct \,\, sum).  
\end{eqnarray*}     
\end{lemma}

\noindent
{\it Proof:} By Lemma \ref{thm:refineU_i} and since $U_0, U_1, \ldots, U_d$ is a decomposition of $V$,\begin{eqnarray*}
V=\sum_{i=0}^{d}\sum_{j=0}^{min(i,d-i)}R^{i-j}H_j \qquad (direct \,\, sum).  
\end{eqnarray*}
In this sum we interchange the order of summation and find 
\begin{eqnarray*}
V=\sum_{i=0}^{d/2}\sum_{j=0}^{d-2i}R^{j}H_i \qquad (direct \,\, sum).  
\end{eqnarray*}
The result now follows by Lemma \ref{thm:scriptDR^iH_i}.
\hfill $\Box$ \\ 

\begin{lemma}
\label{thm:wsumv}
With reference to Definition \ref{def:W}, 
\begin{eqnarray*}
\label{wdirectsumv1}
V=\sum_{i=0}^{d}W_i.
\end{eqnarray*}
\end{lemma}

\noindent
{\it Proof:} Define $V'=\sum_{i=0}^{d}W_i$.  We show $V=V'$.  By construction $V' \subseteq V$.  We now show $V \subseteq V'$.  For  $0 \leq i \leq d$ we set $j=d-i$ in Lemma \ref{thm:AK^{-1}Wij}(i) and find $(A-q^{d-2i}I)W_i \subseteq W_{i+1}$.  Using this we find $AV' \subseteq V'$.  By this and Definition \ref{def:scriptD} we find $\mathcal{D}V' \subseteq V'$.  By Definition \ref{def:W} and Lemma \ref{thm:HVU} we find $H_j \subseteq W_j$ for $0 \leq j \leq d/2$.  Therefore $H_j \subseteq V'$ for $0 \leq j \leq d/2$.  By these comments we find $\mathcal{D}H_j \subseteq V'$ for $0 \leq j \leq d/2$.  Now $V \subseteq V'$ in view of Lemma \ref{thm:VscriptDH_i}.  We have now shown $V=V'$.  
\hfill $\Box$ \\

\begin{corollary}
\label{thm:WVUfinal}
With reference to Assumption \ref{ass}, Definition \ref{def:V}, and Definition \ref{def:W}, the following (i)--(iii) hold.
\begin{enumerate}
\item $W_0+ \cdots +W_i=U_0+ \cdots +U_i$, \qquad $0 \leq i \leq d$,
\item $W_i+ \cdots +W_d=V_0+ \cdots +V_{d-i}$, \qquad $0 \leq i \leq d$,
\item $\dim(W_i)=\rho_i$, \qquad $0 \leq i \leq d$.  
\end{enumerate}
\end{corollary}

\noindent
{\it Proof:} (i) Let $i$ be given.  Define $\Delta = W_0+ \cdots +W_i$ and $\Gamma = U_0+ \cdots +U_i$.  We show $\Delta = \Gamma$.  By Definition \ref{def:W} we have $\Delta \subseteq \Gamma$.  Thus it suffices to show dim$(\Delta)$ $=$ dim$(\Gamma)$.  By construction dim$(\Delta)$ $\leq$ dim$(\Gamma)$.  Suppose, towards a contradiction, that dim$(\Delta)$ $<$ dim$(\Gamma)$.  Using Definition \ref{def:rho}, Lemma \ref{thm:directsum} and since $U_0, U_1, \ldots, U_d$ is a decomposition of $V$ we find
\begin{eqnarray}
\label{stop1}
\sum_{h=0}^{i} \mbox{dim}(W_h) < \sum_{h=0}^{i} \rho_h.
\end{eqnarray} 
By Definition \ref{def:W} we have $\sum_{h=i+1}^dW_h \subseteq \sum_{h=0}^{d-i-1}V_h$.  Using Definition \ref{def:rho}, Lemma \ref{thm:rhofacts}(iii), Lemma \ref{thm:directsum} and since $V_0, V_1, \ldots, V_d$ is a decomposition of $V$ we find  
\begin{eqnarray}
\label{stop2}
\sum_{h=i+1}^{d} \mbox{dim}(W_{h}) \leq \sum_{h=i+1}^{d}\rho_h.
\end{eqnarray} 
By Lemma \ref{thm:directsum} and Lemma \ref{thm:wsumv} we find 
\begin{eqnarray}
\label{stop3}
\mbox{dim}(V) = \sum_{h=0}^{d} \mbox{dim}(W_h).
\end{eqnarray}
Adding (\ref{stop1})--(\ref{stop3}) we find $\mbox{dim}(V)$ $<$ $\sum_{h=0}^{d}\rho_h$.  This contradicts Lemma \ref{thm:rhofacts}(ii) and the result follows.  \\
(ii) Similar to (i). \\
(iii) By (i), Lemma \ref{thm:directsum}, and since $U_0, U_1, \ldots, U_d$ is a decomposition of $V$ we find $\sum_{h=0}^{i} \mbox{dim}(W_h)=\sum_{h=0}^{i} \rho_h$ for $0 \leq i \leq d$.  The result follows.
\hfill $\Box$ \\

\begin{corollary}
\label{thm:wnot0}
With reference to Definition \ref{def:W},
\begin{eqnarray*}
W_i \neq 0, \qquad 0 \leq i \leq d.
\end{eqnarray*}
\end{corollary}

\noindent
{\it Proof:} Immediate from Lemma \ref{thm:rhofacts}(i) and Corollary \ref{thm:WVUfinal}(iii).
\hfill $\Box$ \\

\noindent
Combining Lemma \ref{thm:directsum}, Lemma \ref{thm:wsumv}, and Corollary \ref{thm:wnot0} we obtain Theorem \ref{thm:Wfinal}.  

\section{Interchanging $R$ and $L$}

\noindent
In this section we use Note \ref{note:assumptions} to obtain results that are analogous to the results in Sections 5 and 6.

\begin{definition}
\label{def:A^*}
\rm
With reference to Assumption \ref{ass}, let $A^*:V\to V$ denote the following linear transformation:
\begin{eqnarray*} 
 A^*=K^{-1}+L.
\end{eqnarray*}
\end{definition}

\begin{lemma}
\label{thm:A^*diag}
With reference to Definition \ref{def:A^*} and Assumption \ref{ass}, the following holds. 
$A^*$ is diagonalizable on $V$ and the set of distinct eigenvalues of $A^*$ is $\{\,q^{2i-d} \, |\, 0 \leq i \leq d \, \}$.  Moreover, for $0 \leq i \leq d$, the dimension of the eigenspace for $A^*$ associated with $q^{d-2i}$ is equal to the dimension of $U_i$. 
\end{lemma}

\noindent
{\it Proof:} Apply Note \ref{note:assumptions} to Lemma \ref{thm:Adiag}.  
\hfill $\Box$ \\

\begin{definition}
\label{def:V^*}
\rm
With reference to Definition \ref{def:A^*} and Lemma \ref{thm:A^*diag}, for $0 \leq i \leq d$ we let $V^*_i$ denote the eigenspace for $A^*$ with eigenvalue $q^{d-2i}$.  For notational convenience we set $V^*_{-1}:=0,V^*_{d+1}:=0$.  We observe that $V^*_0, V^*_1, \ldots, V^*_d$ is a decomposition of $V$.   
\end{definition}

\begin{definition}
\label{def:W^*}
\rm  
With reference to Assumption \ref{ass} and Definition \ref{def:V^*}, we define
\begin{eqnarray*}
W^*_i=(U_{i}+ \cdots +U_d) \cap (V^*_{d-i}+ \cdots +V^*_{d}),  \qquad 0 \leq i \leq d.
\end{eqnarray*}
For notational convenience we set $W^*_{-1}:=0, W^*_{d+1}:=0$.
\end{definition}

\begin{theorem}
\label{thm:W^*final}
With reference to Definition \ref{def:W^*}, the sequence $W^*_0, W^*_1, \ldots, W^*_d$ is a decomposition of $V$.
\end{theorem}

\noindent
{\it Proof:} Apply Note \ref{note:assumptions} to Theorem \ref{thm:Wfinal}.
\hfill $\Box$ \\

\section{The linear transformations $B$ and $B^*$}

In this section we introduce the linear transformations $B$, $B^*$ and present a number of relations involving $A$, $A^*$, $B$, $B^*$, $K$, $K^{-1}$.
\begin{definition} 
\label{def:BB^*}
\rm
With reference to Definition \ref{def:W} and Definition \ref{def:W^*}, we define the following linear transformations. 
\begin{enumerate}
\item Let $B:V\to V$ be the linear transformation such that for $0 \leq i \leq d$, $W_i$  is an eigenspace for $B$ with eigenvalue $q^{2i-d}$. 
\item Let $B^*:V\to V$ be the linear transformation such that for $0 \leq i \leq d$, $W_i^*$  is an eigenspace for $B^*$ with eigenvalue $q^{d-2i}$. 
\end{enumerate}
\end{definition}

\begin{lemma}
\label{thm:AA^*BB^*}
{\rm \cite[Lemma 7.2]{Benkart}}
With reference to Definition \ref{def:A}, Definition \ref{def:A^*}, and Definition \ref{def:BB^*},
\begin{eqnarray}
\frac{qAB-q^{-1}BA}{q-q^{-1}} = I,
\label{AA^*BB^*1} \\
\frac{qA^*B^*-q^{-1}B^*A^*}{q-q^{-1}} = I,
\label{AA^*BB^*2} \\
\frac{qBA^*-q^{-1}A^*B}{q-q^{-1}} = I,
\label{AA^*BB^*3} \\
\frac{qB^*A-q^{-1}AB^*}{q-q^{-1}} = I.
\label{AA^*BB^*4}
\end{eqnarray}
\end{lemma}

\begin{lemma}
\label{thm:BB^*K}
{\rm \cite[Lemma 9.1]{Benkart}}
With reference to Assumption \ref{ass} and Definition \ref{def:BB^*},
\begin{eqnarray}
\frac{qBK^{-1}-q^{-1}K^{-1}B}{q-q^{-1}} = I,
\label{BB^*K1} \\
\frac{qB^*K-q^{-1}KB^*}{q-q^{-1}} = I.
\label{BB^*K2}
\end{eqnarray}
\end{lemma}

\begin{lemma}
\label{thm:qserreBB^*}
{\rm \cite[Lemma 10.1]{Benkart}}
With reference to Definition \ref{def:BB^*}, the following (i),(ii) hold.
\begin{enumerate}
\item $B^3B^*-\lbrack 3 \rbrack B^{2}B^*B+\lbrack 3 \rbrack BB^*B^2-B^*B^3=0$,
\item $B^{*3}B-\lbrack 3 \rbrack B^{*2}BB^*+\lbrack 3 \rbrack B^*BB^{*2}-BB^{*3}=0$.
\end{enumerate}
\end{lemma}

\section{The proof of Theorem \ref{thm:main} (existence)}

In this section we prove the existence part of Theorem \ref{thm:main}.  

\begin{definition}
\label{def:rl}
\rm
\noindent
With reference to Assumption \ref{ass} and Definition \ref{def:BB^*}, let $r:V \to V$ and $l:V \to V$ denote the following linear transformations:
\begin{eqnarray*}
r=\frac{I-KB^*}{q(q-q^{-1})^2}, \qquad l=\frac{I-K^{-1}B}{q(q-q^{-1})^2}.
\end{eqnarray*}
\end{definition}

\begin{lemma}
\label{thm:BB^*rl}
With reference to Definition \ref{def:rl}, the following (i),(ii) hold.
\begin{enumerate}
\item $B=K-q(q-q^{-1})^2Kl$,
\item $B^*=K^{-1}-q(q-q^{-1})^2K^{-1}r$.
\end{enumerate}
\end{lemma}

\noindent
{\it Proof:} Immediate from Definition \ref{def:rl}.
\hfill $\Box$ \\

\begin{theorem}
\label{thm:allrelations}
With reference to Assumption \ref{ass} and Definition \ref{def:rl}, the following (i)--(ix) hold.
\begin{enumerate}
\item $KK^{-1}=K^{-1}K=I$, 
\item $KR=q^2RK$, \qquad $KL=q^{-2}LK$, 
\item $Kr=q^2rK$, \quad  $Kl=q^{-2}lK$,
\item $rR=Rr$, \qquad  $lL=Ll$,
\item $lR-Rl=\frac{K^{-1}-K}{q-q^{-1}}$, \qquad $rL-Lr=\frac{K-K^{-1}}{q-q^{-1}}$, 
\item $R^3L-\lbrack 3 \rbrack R^2LR+\lbrack 3 \rbrack RLR^2-LR^3=0$,
\item $L^3R-\lbrack 3 \rbrack L^2RL+\lbrack 3 \rbrack LRL^2-RL^3=0$,
\item $r^3l-\lbrack 3 \rbrack r^2lr+\lbrack 3 \rbrack rlr^2-lr^3=0$,
\item $l^3r-\lbrack 3 \rbrack l^2rl+\lbrack 3 \rbrack lrl^2-rl^3=0$.
\end{enumerate}
\end{theorem}

\noindent
{\it Proof:} (i) Immediate from Note \ref{def:K^-1}. \\
(ii) These equations hold by Lemma \ref{thm:KLR}. \\
(iii) Evaluate Lemma \ref{thm:BB^*K} using Lemma \ref{thm:BB^*rl}.  \\
(iv) Evaluate  (\ref{AA^*BB^*3}),(\ref{AA^*BB^*4}) using Definition \ref{def:A}, Definition \ref{def:A^*}, and Lemma \ref{thm:BB^*rl} and simplify the result using (ii),(iii) above.  \\
(v) Evaluate  (\ref{AA^*BB^*1}),(\ref{AA^*BB^*2}) using Definition \ref{def:A}, Definition \ref{def:A^*}, and Lemma \ref{thm:BB^*rl} and simplify the result using (ii),(iii) above.  \\
(vi),(vii) These relations hold by Assumption \ref{ass}(v),(vi). \\
(viii), (ix) Evaluate Lemma \ref{thm:qserreBB^*} using Lemma \ref{thm:BB^*rl} and simplify the result using (iii) above.
\hfill $\Box$ \\

\begin{theorem}
\label{thm:existence}
With reference to Assumption \ref{ass} and Definition \ref{def:rl}, $V$ supports a $\uqhat$-module structure for which the Chevalley generators act as follows: \\\\
\centerline{
\begin{tabular}[t]{c|cccccccc}
        {\rm generator}  
        &  $e_0^-$  
         & $e_1^-$
         & $e_0^+$  
         & $e_1^+$  
         & $K_0$  
         & $K_1$  
         & $K_0^{-1}$  
         & $K_1^{-1}$  
        \\
        \hline 
{\rm action on $V$} 
& $L$ & $R$ & $r$ & $l$ &
$K$ & $K^{-1}$ &
$K^{-1}$ 
& $K$ 
\end{tabular}}
\end{theorem}

\noindent
{\it Proof:} To see that the above action on $V$ determines a $\uqhat$-module, compare the equations in Theorem \ref{thm:allrelations} with the defining relations for $\uqhat$ in Definition \ref{def:uq}.  
\hfill $\Box$ \\

\noindent 
{\bf Proof of Theorem \ref{thm:main} (existence)}:
The existence part of Theorem \ref{thm:main} is immediate from Theorem \ref{thm:existence}.  
\hfill $\Box$

\medskip
\noindent
Note that the $\uqhat$-module structure given by Theorem \ref{thm:main} is not necessarily irreducible.

\section{The proof of Theorem \ref{thm:main} (uniqueness)} 
 
In this section we prove the uniqueness part of Theorem \ref{thm:main}.  

\medskip
\noindent
The quantum algebra $U_q(sl_2)$ and its finite dimensional modules will be useful in proving uniqueness.  We begin by recalling the definition of $U_q(sl_2)$. 

\begin{definition}
\label{def:uqsl2}
\rm
\cite[Definition 1.1]{Jantzen}
The quantum algebra $U_q(sl_2)$ is the unital associative $\K$-algebra with generators $k,k^{-1}, e, f$ which satisfy the following relations:
\begin{eqnarray*}
kk^{-1}=k^{-1}k=1, \\
ke=q^2ek, \\
kf=q^{-2}fk, \\
ef-fe=\frac{k-k^{-1}}{q-q^{-1}}.
\end{eqnarray*}

\end{definition}

\noindent
We now recall the finite dimensional irreducible $U_q(sl_2)$-modules.   

\begin{lemma}
\label{thm:uqsl2modules}
{\rm \cite[Theorem 2.6]{Jantzen}}
With reference to Definition \ref{def:uqsl2}, there exists a family 
\begin{eqnarray*}
V_{\epsilon, d}, \qquad \epsilon \in \{-1,1\}, \qquad  d=0,1,2,\ldots 
\end{eqnarray*}
of finite dimensional irreducible $U_q(sl_2)$-modules with the following properties.  The module $V_{\epsilon, d}$ has a basis $u_0,u_1, \ldots, u_d$ satisfying:  
\begin{eqnarray}
ku_i = \epsilon q^{d-2i} u_i, \qquad 0 \leq i \leq d, \qquad
\label{ku_i}
\\
fu_i = \lbrack i+1 \rbrack u_{i+1}, \qquad 0 \leq i \leq d-1, \qquad fu_d=0,
\label{fu_i}
\\
eu_i = \epsilon \lbrack d-i+1 \rbrack u_{i-1}, \qquad 1 \leq i \leq d, \qquad eu_0=0.
\label{eu_i} 
\end{eqnarray}
Moreover, every finite dimensional irreducible $U_q(sl_2)$-module is isomorphic to exactly one of the modules $V_{\epsilon, d}$.   
\end{lemma}

\begin{remark}
\label{char2}
\rm
If the characteristic Char$(\K)=2$ then in Lemma \ref{thm:uqsl2modules} we view $\{-1,1\}$ as having a single element.
\end{remark}
 
\begin{lemma}
\label{thm:charknot2}
{\rm \cite[Proposition 2.3]{Jantzen}}
Let $V$ denote a finite dimensional $U_q(sl_2)$-module.  If Char$(\K)\neq 2$ then the action of $k$ on $V$ is diagonalizable.
\end{lemma}

\begin{remark}
\label{counterexample}
\rm
\cite[p.~19]{Jantzen}
Assume Char$(\K)= 2$.  We display a finite dimensional $U_q(sl_2)$-module on which the action of $k$ is not diagonalizable.  Let $k,k^{-1},e,f$ denote the generators for $U_q(sl_2)$ as in Definition \ref{def:uqsl2}.  Let $k$ and $k^{-1}$ act on the vector space $\K^{2}$ as $\begin{pmatrix}1 & 1 \\ 0 & 1 \end{pmatrix}$ and let $e$ and $f$ act on $\K^{2}$ as $0$.  Then $\K^{2}$ is a finite dimensional $U_q(sl_2)$-module and the action of $k$ on $\K^{2}$ is not diagonalizable.
\end{remark}

\begin{lemma}
\label{thm:compred}
{\rm \cite[Theorem 2.9]{Jantzen}}
Let $V$ denote a finite dimensional $U_q(sl_2)$-module.  If the action of $k$ on $V$ is diagonalizable then $V$ is the direct sum of irreducible $U_q(sl_2)$-submodules. 
\end{lemma}

\begin{lemma}
\label{thm:ef0}
Let $V$ be a finite dimensional $U_q(sl_2)$-module.  Assume the action of $k$ on $V$ is diagonalizable.  For $\epsilon \in \{1,-1\}$ and for an integer $d \geq 0$ let $v \in V$ denote an eigenvector for $k$ with eigenvalue $\epsilon q^{d}$.  Then $ev=0$ if and only if $f^{d+1}v=0$ . 
\end{lemma} 

\noindent
{\it Proof:}  Immediate from Lemma \ref{thm:uqsl2modules} and Lemma \ref{thm:compred}. 
\hfill $\Box$ \\

\begin{lemma}
\label{thm:ee'}
Let $V$ denote a finite dimensional vector space over $\K$.  Suppose there are two $U_q(sl_2)$-module structures on $V$.  Assume the actions of $k$ on $V$ given by the two module structures agree and are diagonalizable.  Assume the actions of $f$ on $V$ given by the two module structures agree.  Then the actions of $e$ on $V$ given by the two module structures agree.  
\end{lemma}

\noindent
{\it Proof:}  Let $E:V \to V$ (resp. $E':V \to V$) denote the action of $e$ on $V$ given by the first (resp. second) module structure.  We show $(E-E')V=0$.  Using Lemma \ref{thm:compred} and refering to the first module structure we find $V$ is the direct sum of irreducible $U_q(sl_2)$-submodules.  Let $W$ be one the irreducible submodules in this sum.  It suffices to show $(E-E')W=0$.  By Lemma \ref{thm:uqsl2modules}, there exists a nonnegative integer $d$ and $\epsilon \in \{1,-1\}$ such that $W$ is isomorphic to $V_{\epsilon, d}$.  Therefore, the eigenvalues for $k$ on $W$ are $\epsilon q^{d-2i}$ $(0 \leq i \leq d)$, and dim$(W)=d+1$.  Let $u \in W$ be an eigenvector for $k$ with eigenvalue $\epsilon q^d$.  By Lemma \ref{thm:uqsl2modules}, the vectors $u, fu, \ldots, f^{d}u$ are a basis for $W$.  We show $(E-E')f^{i}u=0$ for $0 \leq i \leq d$.  First assume $i=0$.  Using Lemma \ref{thm:uqsl2modules} we find $Eu=0$.  Also by Lemma \ref{thm:uqsl2modules} we find $f^{d+1}u=0$ and so $E'u=0$ by Lemma \ref{thm:ef0}.  We have now shown $(E-E')u=0$.  Next let $i \geq 1$.  By induction on $i$ we may assume 
\begin{eqnarray}
\label{@1}
(E-E')f^{i-1}u=0.
\end{eqnarray}
Using Definition \ref{def:uqsl2} and since the actions of $k$ (resp. $f$) on $V$ given by the two module structures agree we find $Ef-fE=E'f-fE'$. Hence
\begin{eqnarray}
\label{@2}
f(E-E')=(E-E')f.  
\end{eqnarray}
Applying $f$ to (\ref{@1}) and using (\ref{@2}) we find $(E-E')f^{i}u=0$.  This shows $(E-E')W=0$ and the result follows.    
\hfill $\Box$ \\

\noindent
The following lemma relates $\uqhat$-modules to $U_q(sl_2)$-modules. 

\begin{lemma}
\label{thm:uquqhat}
Let $V$ denote a finite dimensional $\uqhat$-module.  Then for $i \in \{0,1\}$, $V$ supports a $U_q(sl_2)$-module structure such that $K_i-k$, $e_i^{+}-e$, and $e_i^{-}-f$ vanish on $V$.
\end{lemma}

\noindent
{\it Proof:} Immediate from Definition \ref{def:uq} and Definition \ref{def:uqsl2}.  
\hfill $\Box$ \\

\noindent 
{\bf Proof of Theorem \ref{thm:main} (uniqueness)}:
Suppose there exist two $\uqhat$-module structures on $V$ satisfying the conditions of Theorem \ref{thm:main}.  We show that these two module structures agree.  By construction the actions of $e_0^{-}$ (resp. $e_1^{-}, K_0, K_1$) on $V$ given by the two $\uqhat$-module structures agree.  We now show that the actions of $e_0^{+}$ on $V$ given by the two $\uqhat$-module structures agree.  Note that the actions of $K_0$ on $V$ given by the two $\uqhat$-module structures are diagonalizable.  Using Lemma \ref{thm:uquqhat}, $V$ supports two $U_q(sl_2)$-module structures given by the two actions of $K_0, e_0^{+}, e_0^{-}$ on $V$.  Using Lemma \ref{thm:ee'}, we find that the actions of $e_0^{+}$ on $V$ given by the two $\uqhat$-module structures agree.  Similarly, the actions of $e_1^{+}$ on $V$ given by the two $\uqhat$-module structures agree.  We have now shown the two $\uqhat$-module structures of $V$ agree.                
\hfill $\Box$ \\

\section{Which $\uqhat$-modules arise from Theorem \ref{thm:main}?}

\noindent
Theorem \ref{thm:main} gives a way to construct finite dimensional $\uqhat$-modules. Not all finite dimensional $\uqhat$-modules arise from this construction; in this section we determine which ones do. 

\begin{definition}
\label{def:basic}
\rm
\noindent
Let $V$ denote a finite dimensional $\uqhat$-module.  Let $d$ denote a nonnegative integer.  We say $V$ is {\it basic of diameter $d$} whenever there exists a decomposition $U_0, U_1, \ldots, U_d$ of $V$ and linear transformations $R:V\to V$ and $L:V\to V$ satisfying Assumption \ref{ass} such that the given $\uqhat$-module structure on $V$ agrees with the $\uqhat$-module structure on $V$ given by Theorem \ref{thm:main}.  
\end{definition}

\noindent
Our goal for this section is to determine which $\uqhat$-modules are basic.  We begin with a lemma.

\begin{lemma}
\label{thm:K_0diag}
Let $V$ denote a finite dimensional $\uqhat$-module.  If Char$(\K)\neq 2$ then the actions of $K_0$ and $K_1$ on $V$ are diagonalizable.
\end{lemma}

\noindent
{\it Proof:}  For $i \in \{0,1\}$, view $V$ as a $U_q(sl_2)$-module under the action of $K_i, e_i^{+}, e_i^{-}$ as in Lemma \ref{thm:uquqhat}.  The result now follows immediately by Lemma \ref{thm:charknot2}. 
\hfill $\Box$ \\

\begin{remark}
\label{counterexample2}
\rm
Assume Char$(\K)=2$.  We display a finite dimensional $\uqhat$-module on which the actions of $K_0$ and $K_1$ are not diagonalizable.  Let $e_i^{\pm}, K_i^{\pm1}$, $i \in \{0,1\}$, denote the Chevalley generators for $\uqhat$ as in Definition \ref{def:uq}.  Let $K_0^{\pm1}$ and $K_1^{\pm1}$ act on the vector space $\K^{2}$ as $\begin{pmatrix}1 & 1 \\ 0 & 1 \end{pmatrix}$ and let $e_0^{\pm}$ and $e_1^{\pm}$ act on $\K^{2}$ as $0$.  Then $\K^{2}$ is a finite dimensional $\uqhat$-module and the actions of $K_0$ and $K_1$ on $\K^{2}$ are not diagonalizable.
\end{remark}

\begin{theorem}
\label{thm:basiccrit}
Let $d$ denote a nonnegative integer and let $V$ denote a finite dimensional $\uqhat$-module.  With reference to Definition \ref{def:basic}, the following are equivalent.  
\begin{enumerate}
\item $V$ is basic of diameter $d$.
\item $(K_0K_1-I)V=0$,  the action of $K_0$ on $V$ is diagonalizable, and the set of distinct eigenvalues for $K_0$ on $V$ is $\{q^{2i-d}\, | \, 0 \leq i \leq d \}$.
\end{enumerate}
\end{theorem}

\noindent
{\it Proof:} (i) $\Rightarrow$ (ii): Let $U_0, U_1, \ldots, U_d$ be the decomposition of $V$ from Definition \ref{def:basic}.  Let $K:V\to V$ be the linear transformation that corresponds to this decomposition as in Definition \ref{def:K}.  By Definition \ref{def:basic} and Theorem \ref{thm:main}, we find $K-K_0$ and $K^{-1}-K_1$ vanish on $V$.  The result follows. \\
(ii) $\Rightarrow$ (i):  For $0 \leq i \leq d$ let $U_i \subseteq V$ be the eigenspace for the action of $K_0$ on $V$ with eigenvalue $q^{2i-d}$.  Note that the sequence $U_0, U_1, \ldots, U_d$ is a decomposition of $V$.  Define linear transformations $R:V\to V$ and $L:V\to V$ by $(R-e_{1}^{-})V=0$ and $(L-e_{0}^{-})V=0$.  We now check that $R$ and $L$ satisfy Assumption \ref{ass}(i)--(vi).  Using (\ref{eq:buq3}) and (\ref{eq:buq4}), it is routine to check that $R$ and $L$ satisfy Assumption \ref{ass}(i),(ii).  Next we verify that $R$ satisfies Assumption \ref{ass}(iii).  View $V$ as a $U_q(sl_2)$-module under the action of $K_1, e_1^+, e_1^-$ as in Lemma \ref{thm:uquqhat}.  Since $(K_0K_1-I)V=0$ and the action of $K_0$ on $V$ is diagonalizable we find that the action of $K_1$ on $V$ is diagonalizable.  Thus, by Lemma \ref{thm:compred}, $V$ is a direct sum of irreducible $U_q(sl_2)$-submodules.  Let $W$ be one of the irreducible submodules in this sum.  Note that $RW \subseteq W$.  Using Lemma \ref{thm:uqsl2modules}, we find that for $0 \leq i \leq d/2$ the restriction $R^{d-2i}|_{W \cap U_{i}}:W \cap U_{i}\to W \cap U_{d-i}$ is a bijection.  It follows that for $0 \leq i \leq d/2$ the restriction $R^{d-2i}|_{U_{i}}:U_{i}\to U_{d-i}$ is a bijection.  We have now shown $R$ satisfies Assumption \ref{ass}(iii).  The proof that $L$ satisfies Assumption \ref{ass}(iv) is similar.  By (\ref{eq:buq7}) we find $R$ and $L$ satisfy Assumptions \ref{ass}(v),(vi).  We have now shown that $R$ and $L$ satisfy Assumption \ref{ass}(i)--(vi) and so Theorem \ref{thm:main} applies.  It remains to show that the given $\uqhat$-module structure on $V$ agrees with the $\uqhat$-module structure on $V$ given by Theorem \ref{thm:main}.  Using the uniqueness statement in Theorem \ref{thm:main} it suffices to show each of $R-e_{1}^{-}$, $L-e_{0}^{-}$, $K-K_0$, $K^{-1}-K_1$ vanish on $V$, where $K:V\to V$ is the linear transformation that corresponds to the decomposition $U_0, U_1, \ldots, U_d$ as in Definition \ref{def:K}.  The first two equations were mentioned earlier.  Using Definition \ref{def:K} we find $(K-K_0)V=0$.  Since $(K_0K_1-I)V=0$ we find $(K^{-1}-K_1)V=0$.  We have now shown that the given $\uqhat$-module structure on $V$ agrees with the $\uqhat$-module structure on $V$ given by Theorem \ref{thm:main}.\hfill $\Box$ \\

\section{The relationship between general $\uqhat$-modules and basic $\uqhat$-modules}

\noindent
Throughout this section $V$ will denote a nonzero finite dimensional $\uqhat$-module (not necessarily irreducible) on which the actions of $K_0$ and $K_1$ are diagonalizable (see Lemma \ref{thm:K_0diag}).   

\medskip
\noindent
In this section we will show, roughly speaking, that $V$ is made up of basic $\uqhat$-modules.  We will use the following definition. 

\begin{definition}
\label{def:8pieces}
\rm
\noindent
For $\epsilon_0, \epsilon_1 \in \{1,-1\}$ we define 
\begin{eqnarray*}
V_{even}^{(\epsilon_0, \epsilon_1)}=\mbox{span}\{\,v \in V\, |\, K_0v=\epsilon_0q^{i}v, \, K_1v=\epsilon_1q^{-i}v, \, i \in \Z, \, i \, even \,\}, 
\end{eqnarray*}
\begin{eqnarray*}
V_{odd}^{(\epsilon_0, \epsilon_1)}=\mbox{span}\{\,v \in V\, |\, K_0v=\epsilon_0q^{i}v, \, K_1v=\epsilon_1q^{-i}v, \, i \in \Z, \, i \, odd \,\}. 
\end{eqnarray*}
\end{definition}

\begin{theorem}
\label{thm:directsum8pieces}
With reference to Definition \ref{def:8pieces},   
\begin{eqnarray}
\label{blob}
V=\sum_{(\epsilon_0,\epsilon_1)}\sum_{\sigma}V_{\sigma}^{(\epsilon_0,\epsilon_1)} \qquad (direct \,\, sum \,\, of \,\, \uqhat-modules),
\end{eqnarray}
where the first sum is over all ordered pairs $(\epsilon_0,\epsilon_1)$ with $\epsilon_0,\epsilon_1 \in \{1,-1\}$, and the second sum is over all $\sigma \in \{even, odd\}$.  
\end{theorem}  

\noindent
{\it Proof:}  Using (\ref{eq:buq2}) and Definition \ref{def:8pieces} we find $V_{even}^{(\epsilon_0, \epsilon_1)}$ and $V_{odd}^{(\epsilon_0, \epsilon_1)}$ are invariant under the action of $K_0^{\pm1}$ and $K_1^{\pm1}$.  Using (\ref{eq:buq3}) and (\ref{eq:buq4}) we find $V_{even}^{(\epsilon_0, \epsilon_1)}$ and $V_{odd}^{(\epsilon_0, \epsilon_1)}$ are invariant under the action of $e_{0}^{\pm}, e_{1}^{\pm}$.  Thus the subspaces on the right hand side of (\ref{blob}) are $\uqhat$-submodules of $V$.  We now show the sum on the right hand side of (\ref{blob}) equals $V$.  Recall that the actions of $K_0$ and $K_1$ on $V$ are both diagonalizable.  Using this and (\ref{eq:buq2}) we find that the actions of $K_0$ and $K_1$ on $V$ are simultaneously diagonalizable.  Thus $V$ is the direct sum of common eigenspaces for the actions of $K_0$ and $K_1$ on $V$.  It remains to show that any common eigenvector for the actions of $K_0$ and $K_1$ on $V$ is in one of the subspaces on the right hand side of (\ref{blob}).  Let $v$ denote a common eigenvector for the action of $K_0$ and $K_1$ on $V$.  By construction there exists $\alpha, \beta \in \K$ such that $K_0v=\alpha v$ and $K_1v=\beta v$.  Using Lemma \ref{thm:uquqhat}, Lemma \ref{thm:compred}, and Lemma \ref{thm:uqsl2modules} we find that there exists an $\epsilon_0 \in \{1,-1\}$ and an integer $i$ such that $\alpha =\epsilon_0q^{i}$.  For every $m \in \Z$ define $T_m:=\{\,x \in V\, |\, K_0x=\epsilon_0 q^{i+2m}x$ and $K_1x=\beta q^{-2m}x \,\}$, and define $T:=\sum_{m \in \Z}T_m$.  Observe $K_0K_1-\epsilon_0 q^i \beta I$ vanishes on $T$.  Using (\ref{eq:buq3}) and (\ref{eq:buq4}) we find $T$ is a $\uqhat$-module.  Also, $0 \neq v \in T$ and so $T$ is not the zero module.  Let $W$ denote an irreducible $\uqhat$-module contained in $T$.  By \cite[Proposition 3.2]{cp}, there exists an $\epsilon \in \{1,-1\}$ such that $K_0K_1-\epsilon I$ vanishes on $W$.  Also, $K_0K_1-\epsilon_0 q^i \beta I$ vanishes on $W$.  So we find $\epsilon=\epsilon_0 q^{i} \beta$.  Define $\epsilon_1=\epsilon \epsilon_0^{-1}$.  Then $\epsilon_1 \in \{1,-1\}$ and $\beta=\epsilon_1 q^{-i}$.  We have now shown $K_0v=\epsilon_0 q^iv$ and $K_1v=\epsilon_1 q^{-i}v$.  Therefore $v \in V_{even}^{(\epsilon_0, \epsilon_1)}$ if $i$ is even or $v \in V_{odd}^{(\epsilon_0, \epsilon_1)}$ if $i$ is odd.  This shows that the sum on the right hand side of (\ref{blob}) equals $V$.  By Definition \ref{def:8pieces} the sum in (\ref{blob}) is direct.
\hfill $\Box$ \\

\begin{lemma}
\label{thm:onepiecebasiceven}
With reference to Definition \ref{def:basic}, Definition \ref{def:8pieces}, and Theorem \ref{thm:directsum8pieces}, the following are equivalent.
\begin{enumerate}
\item $V=V_{even}^{(1,1)}$.
\item $V$ is basic of even diameter.
\item The spaces $V_{even}^{(-1,1)}$, $V_{even}^{(1,-1)}$, $V_{even}^{(-1,-1)}$, $V_{odd}^{(1,1)}$, $V_{odd}^{(-1,1)}$, $V_{odd}^{(1,-1)}$, $V_{odd}^{(-1,-1)}$ are all zero. 
\end{enumerate}
\end{lemma}

\noindent
{\it Proof:}  (i)$\Rightarrow$(ii):  We use Theorem \ref{thm:basiccrit}.  By Definition \ref{def:8pieces}, we find $K_0K_1-I$ vanishes on $V$.  Recall that the action of $K_0$ on $V$ is diagonalizable.  Using Lemma \ref{thm:uquqhat}, Lemma \ref{thm:compred} and Lemma \ref{thm:uqsl2modules} we find that there exists a nonnegative integer $d$ such that the set of distinct eigenvalues for the action of $K_0$ on $V$ is $\{q^{2i-d}\, | \, 0 \leq i \leq d \}$.  So by Theorem \ref{thm:basiccrit} $V$ is basic of diameter $d$. By Definition \ref{def:8pieces} $d$ is even.  \\
(ii)$\Rightarrow$(i):  Immediate from Theorem \ref{thm:basiccrit} and Definition \ref{def:8pieces}. \\
(i)$\Leftrightarrow$(iii):  Immediate from Theorem \ref{thm:directsum8pieces}.  
\hfill $\Box$ \\

\begin{lemma}
\label{thm:onepiecebasicodd}
With reference to Definition \ref{def:basic}, Definition \ref{def:8pieces}, and Theorem \ref{thm:directsum8pieces}, the following are equivalent.
\begin{enumerate}
\item $V=V_{odd}^{(1,1)}$.
\item $V$ is basic of odd diameter.
\item The spaces $V_{odd}^{(-1,1)}$, $V_{odd}^{(1,-1)}$, $V_{odd}^{(-1,-1)}$, $V_{even}^{(1,1)}$, $V_{even}^{(-1,1)}$, $V_{even}^{(1,-1)}$, $V_{even}^{(-1,-1)}$ are all zero. 
\end{enumerate}
\end{lemma}

\noindent
{\it Proof:}  Similar to the proof of Lemma \ref{thm:onepiecebasiceven}.
\hfill $\Box$ \\

\noindent 
Refering to (\ref{blob}), even though the six submodules $V_{even}^{(-1,1)}$, $V_{even}^{(1,-1)}$, $V_{even}^{(-1,-1)}$, $V_{odd}^{(-1,1)}$, $V_{odd}^{(1,-1)}$, $V_{odd}^{(-1,-1)}$ are not basic they become basic after a routine normalization.  This is explained in the following lemma, definition, and remark.  

\begin{lemma}
\label{thm:uqauto}
{\rm \cite[Prop. 3.3]{cp}} For any choice of scalars $\epsilon_0, \epsilon_1$ from $\{1,-1\}$, there exists a $\K$-algebra automorphism of $\uqhat$ such that 
\begin{eqnarray*}
K_i \to \epsilon_iK_i, \qquad e_i^+ \to \epsilon_i e_i^+, \qquad e_i^- \to e_i^-,  
\end{eqnarray*}
for $i \in \{0,1\}$.  We refer to the above automorphism as $\tau(\epsilon_0, \epsilon_1)$.  
\end{lemma}

\begin{definition}
\label{def:twisted}
\rm
Let $W$ denote a $\uqhat$-module.  Let $\tau$ be an automorphism of $\uqhat$.  We define a new $\uqhat$-module structure on $W$ as follows.  For $x \in \uqhat$ and $w \in W$ define $x . w \,\, ( {\rm new \,\, action} ) = \tau(x) . w \,\, ( {\rm original \,\, action} )$.  We refer to this new $\uqhat$-module structure as {\it $W$ twisted via $\tau$}.  
\end{definition}

\begin{remark}
\rm
With reference to (\ref{blob}) each submodule $V_{\sigma}^{(\epsilon_0,\epsilon_1)}$ becomes basic upon twisting via $\tau(\epsilon_0, \epsilon_1)$. 
\end{remark}

\section{Acknowledgment} 

\noindent
I would like to express my gratitude to my thesis advisor Paul Terwilliger for introducing me to this subject and for his many useful suggestions.

\bigskip

\noindent 
Darren Funk-Neubauer 
\medskip \hfil\break
\noindent Department of Mathematics \hfil\break
\noindent University of Wisconsin-Madison \hfil\break 
\noindent 480 Lincoln Drive \hfil\break
\noindent Madison, WI 53706-1388 USA 
\medskip

\noindent
email:   {\tt neubauer@math.wisc.edu}


\begin{thebibliography}{10}

\bibitem{Curtin}
H.~Al-Najjar, B.~Curtin, A family of tridiagonal pairs related to the quantum affine algebra $\uqhat$, Electron. J. Linear Algebra 13 (2005) 1--9. 

\bibitem{Beck}
J.~Beck, H.~Nakajima, Crystal bases and two-sided cells of quantum affine algebras, Duke Math. J. 123 (2004) 335--402.

\bibitem{Benkart}
G.~Benkart, P.~Terwilliger, Irreducible modules for the quantum affine algebra $\uqhat$ and its Borel subalgebra, J. Algebra 282 (2004) 172--194.

\bibitem{cp}
V.~Chari, A.~Pressley, Quantum affine algebras, Commun. Math. Phys. 142 (1991) 261--283.

\bibitem{cp2}
V.~Chari, A.~Pressley, Quantum affine algebras and their representations, in: Canad.  Math.  Soc.  Conf.  Proc., vol. 16, Amer. Math. Soc., Providence, RI, 1995, pp. 1659--1678.

\bibitem{Drinfeld1}
V.G.~Drinfeld, Hopf algebras and the quantum Yang-Baxter equation, Soviet Math. Dokl. 32 (1985) 254--258.

\bibitem{Drinfeld2}
V.G.~Drinfeld, A new realization of Yangians and quantized affine algebras, Soviet Math. Dokl. 36 (1988) 212--216.

\bibitem{Ito}
T.~Ito, K.~Tanabe, P.~Terwilliger, Some algebra related to P- and Q-polynomial association schemes, in: Codes and Association Schemes, Piscataway NJ, 1999, in: DIMACS Ser. Discrete Math. Theoret. Comput. Sci., vol. 56, Amer. Math. Soc., Providence, RI, 2001, pp. 167--192.

\bibitem{shape}
T.~Ito, P.~Terwilliger, The shape of a tridiagonal pair, J. Pure Appl. Algebra 188 (2004) 145--160.

\bibitem{Ito.Ter} 
T.~Ito, P.~Terwilliger, Tridiagonal pairs and the quantum affine algebra $\uqhat$, Ramanujan J., in press.  

\bibitem{Jantzen}
J.C.~Jantzen, Lectures on quantum groups, Amer. Math. Soc., Providence RI, 1996.

\bibitem{Jimbo1}
M.~Jimbo, A $q$-difference analogue of $U(\mathfrak{g})$ and the Yang-Baxter equation, Lett. Math. Phys. 10 (1985) 63--69.

\bibitem{Jimbo2}
M.~Jimbo, A $q$-analogue of $U(\mathfrak{gl}(N+1))$, Hecke algebra and the Yang-Baxter equation, Lett. Math. Phys. 11 (1986) 247--252.

\bibitem{Kang}
S-J.~Kang, J-H.~Kwon, Crystal bases of the Fock space representations and string functions, J. Algebra 280 (2004) 313--349.

\bibitem{Mansour}
M.~Mansour, E.H. Zakkari, Fractional spin through quantum affine algebras with vanishing central charge, Internat. J. Theoret. Phys. 43 (2004) 1249--1260.

\bibitem{Misra}
K.~Misra, V.~Williams, Combinatorics of quantum affine Lie algebra representations, in: Kac-Moody Lie algebras and related topics, Contemp. Math., vol. 343, Amer. Math. Soc., Providence, RI, 2004, pp. 167--189.

\bibitem{Thoren}
J.~Thoren, Finite-dimensional modules of quantum affine algebras, J. Algebra 272 (2004) 581--613.
\end{thebibliography}
\end{document}